\documentclass[11pt]{amsart}

\usepackage{fullpage}

\usepackage{url}

\usepackage{amssymb}

\usepackage{cite}

\renewcommand{\a}{\alpha}

\newcommand{\g}{\gamma}
\renewcommand{\d}{\delta}
\newcommand{\D}{\Delta}
\newcommand{\e}{\varepsilon}
\newcommand{\f}{\varphi}

\renewcommand{\S}{\Sigma}
\renewcommand{\k}{\kappa}
\renewcommand{\l}{\lambda}
\newcommand{\z}{\zeta}

\newcommand{\cC}{{\mathcal C}}
\newcommand{\cM}{{\mathcal M}}

\newcommand{\cB}{{\mathcal B}}
\newcommand{\cL}{{\mathcal L}}

\newcommand{\cH}{{\mathcal H}}
\newcommand{\cD}{{\mathcal D}}

\newcommand{\bC}{\mathbb C}
\newcommand{\bR}{\mathbb R}

\newcommand{\bZ}{\mathbb Z}
\newcommand{\bS}{\mathbb S}

\newcommand{\bH}{\mathbb H}

\newcommand{\be}{\begin{equation}}
\newcommand{\ee}{\end{equation}}

\newcommand{\bel}[1]{\begin{equation}\label{#1}}
\newcommand{\beaa}{\begin{eqnarray*}}
\newcommand{\bea}{\begin{eqnarray}}
\newcommand{\beal}[1]{\begin{eqnarray}\label{#1}}
\newcommand{\bean}{\begin{eqnarray}\nonumber}
\newcommand{\beadl}[1]{\begin{deqarr}\label{#1}}
\newcommand{\eeadl}[1]{\arrlabel{#1}\end{deqarr}}
\newcommand{\eeal}[1]{\label{#1}\end{eqnarray}}
\newcommand{\eead}[1]{\end{deqarr}}
\newcommand{\eea}{\end{eqnarray}}
\newcommand{\eeaa}{\end{eqnarray*}}

\renewcommand{\to}{\rightarrow}

\renewcommand{\phi}{\varphi}
\renewcommand{\epsilon}{\varepsilon}

\newcommand{\<}{\langle}
\renewcommand{\>}{\rangle}

\newcommand{\dm}{{\partial M}}

\newcommand{\w}{\widetilde}

\theoremstyle{plain}
\newtheorem{theorem}{Theorem}[section]

\newtheorem{remark}[theorem]{Remark}

\newtheorem{lemma}[theorem]{Lemma}

\newtheorem{proposition}[theorem]{Proposition}

\theoremstyle{definition}

\def\endproof{\qed \medskip}
\def\blacksquare{\hbox to .60em {\vrule width .60em height .60em}}

\numberwithin{equation}{section}

\date{\today}

\begin{document}

\title[ ]{Conformal Immersions of Prescribed Mean Curvature  in $\bR^{3}$}

\author[ ]{Michael T. Anderson}

\address{Department of Mathematics, Stony Brook University, Stony Brook, N.Y.~11794-3651, USA} 
\email{anderson@math.sunysb.edu}
\urladdr{http://www.math.sunysb.edu/$\sim$anderson}

\thanks{Partially supported by NSF grant DMS 0905159 }

\begin{abstract}
We prove the existence of (branched) conformal immersions $F: S^{2} \to \bR^{3}$ with mean curvature 
$H > 0$ arbitrarily prescribed up to a 3-dimensional affine indeterminacy. A similar result is proved 
for the space forms $\bS^{3}$, $\bH^{3}$ and partial results for surfaces of higher genus. 

\end{abstract}

\maketitle 

\setcounter{section}{0}
\setcounter{equation}{0}

\section{Introduction}

   Consider an immersed closed surface 
$$F: S \to \bR^{3},$$ 
in Euclidean space $\bR^{3}$. In terms of metric or Riemannian geometry, the Cauchy data of 
the immersion $F$ consist of the first and second fundamental forms $(\g, A)$ of the surface 
$\S = Im F$. A basic and natural question is what are the possible elliptic data that can be imposed 
on $F$, giving then a well-defined and geometric elliptic PDE problem for the immersion. Roughly speaking, 
elliptic data consist of one-half of the Cauchy data, so a combination of three of the six components 
of $(\g, A)$. (This of course matches the three scalar components of $F$). 

   Clearly one should consider first the cases of prescribing the induced metric $\g$ or the 
second fundamental form $A$. One can view $\g$ as Dirichlet data and $A$ as Neumann 
data for $F$. This becomes clearer when one considers $\S$ as the boundary of compact domain 
$M \subset \bR^{3}$ (or immersion of a domain $M$ into $\bR^{3}$); the pair $(\g, A)$ is then Dirichlet 
and Neumann data for a flat metric on $M$ with $\dm = \S$.

   Prescribing the induced metric $\g$ is the well-known isometric immersion problem: given a 
metric $\g$ on $S$, is there an isometric immersion $F: S \to \bR^{3}$, i.e.~$F^{*}(g_{Eucl}) = \g$? 
This problem has remained notoriously difficult despite much effort. From a large-scale viewpoint, 
relatively little progress has been made on this global problem except in the case of positive 
Gauss curvature $K > 0$, i.e.~the solution of the Weyl problem by Nirenberg \cite{N} and Pogorelov 
\cite{P}. 

  However, prescribing the induced metric $\g$ is not an elliptic problem. This follows easily from 
Gauss' Theorema Egregium,
\be \label{TE}
K = det A,
\ee
where $K$ is the Gauss curvature. Namely, if prescribing a $C^{m,\a}$ metric $\g \in Met^{m,\a}(S)$ 
were elliptic, then by elliptic regularity, $F$ would be a $C^{m+1,\a}$ mapping and hence $A$ would be 
a $C^{m-1,\a}$ form on $S$. Then \eqref{TE} gives $K \in C^{m-1,\a}(S)$. However, the space of $C^{m,\a}$ 
metrics for which the curvature $K$ is in $C^{m-1,\a}$ is of infinite codimension, so that the operator 
$F \to \g = F^{*}(g_{Eucl})$ is not Fredholm. This contradicts ellipticity. Similar reasoning shows 
that prescribing the second fundamental form $A$ is never elliptic. 

  It is perhaps worth noting that isometric immersions can be described by an associated Darboux 
equation, cf.~\cite{HH} for instance, which in case the relation $K > 0$ holds is an elliptic equation of 
Monge-Ampere type. However, there is a loss of one derivative in passing to the Darboux equation, 
so there is no conflict with the lack of ellipticity of the immersion problem itself. 
   
\medskip 

  Returning to the issue of elliptic data for the map $F$, it was shown in \cite{An2} that the data 
$([\g], H)$, where $[\g]$ is the pointwise conformal class of the induced metric $\g$ and 
$H = H_{F}$ is the mean curvature of the immersion $F$, are elliptic data. In fact, to the 
author's knowledge, this is the only known elliptic data, depending only on the Cauchy data 
$(\g, A)$. The data $([\g], H)$ form a nonlinear elliptic system of three equations in the three 
unknowns - the components of the mapping $F$. Note that $[\g]$ involves the first derivatives of 
$F$ while $H$ involves second derivatives. 

\medskip 

  In this paper we study the question of global existence of solutions to this elliptic problem. 

{\it Question}. Given a closed orientable surface $S$ and arbitrary smooth data $([\g], H)$ on $S$, 
does there exist an immersion $F: S \to \bR^{3}$ realizing $([\g], H)$, i.e.
$$([F^{*}(g_{Eucl}), H_{F}]) = ([\g], H).$$

   In this generality, the answer is clearly no. For instance, for any immersion $F$ of a compact 
surface in $\bR^{3}$, one must have $H_{F} > 0$ somewhere. There are in fact further obstructions, 
at least in the case $S = S^{2}$. Namely, the mean curvature $H$ of a conformal 
immersion $(S, [\g]) \to (\bR^{3}, g_{Eucl})$ must satisfy
\be \label{ob}
\int_{S}V(H)dV_{\g} = 0,
\ee
for any conformal vector field $V$ on $S$, cf.~Section 2 for a proof. Thus for instance if $H$ 
is a monotone function of a standard height function $z$ on $S^{2}$, then $H$ cannot be 
realized as the mean curvature of any conformal immersion $S^{2} \to \bR^{3}$. Formally, 
\eqref{ob} gives a 3-dimensional space of conditions on $H$, corresponding to the 3-dimensional 
space of (restrictions of) linear functions on $S^{2}$.

\medskip 

    We mainly focus here on the case $S = S^{2}$; a discussion of the case of higher 
genus is given in Section 4. To state the main result, let $C_{+}^{m-1,\a}$ denote the 
space of $C^{m-1,\a}$ functions $H: S^{2} \to \bR^{+}$, so $H > 0$ everywhere, with 
$m \geq 1, \a \in (0,1)$. Define an equivalence relation on $C_{+}^{m-1,\a}$ by 
\be \label{eq}
[H_{1}] =  [H_{2}]  \Leftrightarrow   H_{2} = H_{1} + \ell,
\ee
where $\ell = a + bx$ and $x$ is the restriction to $\bS^{2}(1) \subset \bR^{3}$ of a linear function on 
$\bR^{3}$ of norm 1 and $a = |b| \geq 0$. Thus $\ell \geq 0$ on $S^{2}$ with $\ell = 0$ at the point $-x \in S^{2}$. 
In particular, the gradients $\nabla \ell$ form the space of conformal vector fields on $S^{2}$. 
Let $\cD_{+}^{m-1,\a} = C_{+}^{m-1,\a} / \sim$ be the quotient space.

\begin{theorem}
  For any pointwise $C^{m,\a}$ conformal class $[\g]$ of metrics on $S^{2}$ and for any equivalence class 
$[H] \in \cD_{+}^{m-1,\a}$, there exists a $C^{m+1,\a}$ branched immersion of $S^{2}$ into $\bR^{3}$ 
with prescribed pointwise conformal class $[\g]$ and prescribed mean curvature class $[H]$. Thus, there 
exists a $C^{m+1,\a}$ smooth branched immersion $F: (S^{2}, \g) \to (\bR^{3}, g_{+1})$ such that 
\be \label{1.2}
[F^{*}(g_{Eucl})] = [\g] \ \ {\rm and} \ \ H(F(x)) = H + \ell,
\ee
for any $H > 0$ and for some affine function $\ell$ as above on $S^{2}$. 
\end{theorem}

  We refer to Section 2 for the precise definition of branched immersion; (it is the usual definition). 
Branched immersions also satisfy the obstruction \eqref{ob}. Generically the data $([\g], H)$ is 
realized by a (regular) immersion; the space of branched immersions corresponds to a ``boundary" 
of the space of immersions which has codimension at least $6$. On the other hand, there is no 
reason to expect that Theorem 1.1 is valid within the space of immersions itself. The presence 
of branch points corresponds to bubbling behavior, commonly arising in geometric PDE problems. 
Such bubble formation can be ruled out for ``low-energy" solutions, corresponding here to 
the mean curvature $H$ close to a constant, but is not to be expected to be ruled out in general; 
see also the discussion of examples in Section 2.

  The literature of studies of surfaces of constant and more generally prescribed mean curvature in 
Euclidean space $\bR^{3}$ and the related space forms $\bS^{3}$ and $\bR^{3}$ is vast. Some of 
the issues discussed here are related to a question posed by Yau \cite{Y}: given a smooth 
function $H: \bR^{3} \to \bR$, is there an immersion $F$ of a surface $\S \to \bR^{3}$ such that 
the mean curvature of $F(\S)$ equals $H$, cf.~\cite{CM}, \cite{GT}, \cite{TW} and further references 
therein for example. Theorem 1.1 deals with a closely related but nevertheless somewhat different 
problem. Also, this paper may be viewed as a continuation of an earlier study in \cite{An2}. 

   Theorem 1.1 is proved in Section 3 by a global degree-theoretic argument, based on the degree of 
proper Fredholm maps between Banach manifolds of Fredholm index 0. In Section 2, we set up the 
basic framework and background, proving in particular that the space of branched immersions 
$F: S^{2} \to \bR^{3}$ is a smooth Banach manifold, for which the map to the target data 
$F \to ([\g], H)$ is smooth and Fredholm, of Fredholm index 0. In Section 3, we prove a basic 
apriori bound on the area of an immersion $F$ in terms of the target data $([\g], H)$. This 
is the essential estimate used to obtain a proper Fredholm map. A computation of the 
degree then follows, based on the Hopf uniqueness theorem for constant mean curvature 
spheres immersed in $\bR^{3}$. 
  
   Theorem 1.1 generalizes naturally to the other space forms $\bS^{3}$ and $\bH^{3}$. This is 
discussed in Section 4, cf.~Theorem 4.1. Also in Section 4, we point out that all of the results 
above carry over to closed orientable surfaces of any genus $g > 0$, except for the last result, namely 
the computation of the degree. This remains an interesting open question.

\section{Preliminary Material}
\setcounter{equation}{0}

 In this section we describe background material and results needed to prove Theorem 1.1.

  Let $Map^{m+1,\a}(S^{2}, \bR^{3})$ denote the space of $C^{m+1,\a}$ maps $F: S^{2} 
\to \bR^{3}$. Throughout the paper we assume $m \geq 1$ and also allow $m = \infty$. This 
is a smooth Banach manifold, in fact Banach space due to the linear structure on $\bR^{3}$; when 
$m = \infty$ one has a Fr\'echet space. Although $Map^{m+1,\a}(S^{2}, \bR^{3})$ is not 
separable, it is separable with respect to a slightly weaker topology, namely the 
$C^{m+1,\a'}$ topology, for any $\a' < \a$, cf.~\cite{W} for further discussion. 
The tangent space at a map $F$ is given by the space of $C^{m+1,\a}$ vector fields 
$X$ along the map $F$. 

  Next, let $Imm^{m+1,\a}(S^{2}, \bR^{3})$ be the space of $C^{m+1,\a}$ immersions $F: S^{2} \to  
\bR^{3}$. This is an open domain in $Map^{m+1,\a}(S^{2}, \bR^{3})$ and so of course is also a 
smooth Banach manifold. Let $\S = Im F$ denote the immersed sphere in $\bR^{3}$. 

   Let $\cC^{m,\a} = \cC^{m,\a}(S^{2})$ be the space of (pointwise) conformal classes $[\g]$ of $C^{m,\a}$ 
metrics on $S^{2}$; recall that two metrics $\g_{1}$ and $\g_{2}$ on $S^{2}$ are conformally equivalent if 
$\g_{2} = \mu^{2}\g_{1}$, for some positive function $\mu$ on $S^{2}$.  Let $C^{m-1,\a}  = C^{m-1,\a}(S^{2})$ 
be the space of $C^{m-1,\a}$ functions on $S^{2}$. An immersion $F \in Imm^{m+1,\a}(S^{2}, \bR^{3})$ 
induces a metric $\g = F^{*}(g_{Eucl})$ on $S^{2}$ and mean curvature function $H = 
H_{F} \in C^{m-1,\a}$, with respect to a choice of normal. Here the normal is chosen 
to be that induced by the outward normal to the sphere tangent to $\S = Im F$ at a point 
where $|F|^{2}$ is maximal.  

   This data may be assembled to a natural map of Banach manifolds  
\be \label{Pi}
\Pi: Imm^{m+1,\a}(S^{2}, \bR^{3}) \to \cC^{m,\a} \times C^{m-1,\a},
\ee
$$\Pi(F) = ([\g], H).$$

   By \cite{An1}, \cite{An2}, the data $([\g], H)$ form an elliptic system of PDEs for the map $F$, 
cf.~also the proof of Proposition 2.2 below. Thus the map $\Pi$ is a Fredholm map, i.e.~the 
linearization $D\Pi$ at any $F$ is a Fredholm linear map. However, the map $\Pi$ has a rather 
large degeneracy, due to the large isometry group $Isom(\bR^{3})$ of $(\bR^{3}, g_{Eucl})$. 
This group acts freely on $Imm \equiv Imm^{m+1,\a}(S^{2}, \bR^{3})$ via 
$(F, \iota) \to \iota \circ F$, corresponding to translation, rotation or reflection of $F$ and 
fixes the target data, i.e.~$\Pi(\iota \circ F) = \Pi(F)$. To remove this degeneracy, we divide 
$Imm$ by this action, and consider only the quotient space $Imm_{b}$ of based immersions. 
There is a global slice to this action, i.e.~an inclusion $Imm_{b} \subset Imm$, given by 
fixing a point $p_{0} \in \bS^{2}(1) \subset \bR^{3}$ (say the north pole), a unit vector 
$e \in T_{p_{0}}(\bS^{2}(1))$ and requiring that $F(p_{0}) = 0$,  $T_{0}(\S) = \bR^{2} 
\subset \bR^{3}$ and with $F_{*}(e) = (a, 0, 0)$ for some $a > 0$. Thus, unless 
mentioned otherwise, throughout the paper we consider
\be \label{Pib}
\Pi: Imm_{b}^{m+1,\a}(S^{2}, \bR^{3}) \to \cC^{m,\a} \times C^{m-1,\a},
\ee
$$\Pi(F) = ([\g], H).$$
By \cite{An1}, the Fredholm index of $\Pi$ equals zero.  (The Fredholm index of $\Pi$ in \eqref{Pi} 
equals $6$, the dimension of $Isom(\bR^{3})$). 

\medskip 

  The basic issue in studying the global properties of the map $\Pi$ in \eqref{Pib} is 
whether $\Pi$ is proper. This is not true per se, due to the non-compactness of the conformal 
group of $S^{2}$. (This is another reason for dividing out by the action of isometries above). 
In fact, by the uniformization theorem, the group of $C^{m+1,\a}$ diffeomorphisms 
of $S^{2}$ acts transitively on the space $\cC^{m,\a}$, with stabilizer the conformal group 
${\mathrm{Conf}}(S^{2})$ of $(S^{2}, [\g])$. This is a non-compact group, diffeomorphic to $\bR^{3}$. 
The group ${\mathrm{Conf}}(S^{2})$ also acts on the space of functions $C^{m-1,\a}$ by pre-composition, 
$(H, \f) \to H\circ \f$, but here the action is proper, except at the constant functions $H = c$. It follows 
that for any point $([\g], c) \in \cC^{m,\a}\times C^{m-1,\a}$, the space $\Pi^{-1}(\cC^{m,\a} \times \{c\})$ 
is non-compact, so that $\Pi$ is not proper. On the other hand, the action of the conformal group 
on $\cC^{m,\a} \times C^{m-1,\a}$ is proper, whenever $H$ is non-constant. 

\medskip 

   There are two ways to deal with this issue and we will in fact use both. First, one may simply remove 
the sets $\Pi^{-1}([\g], c)$ above from the domain of $\Pi$ in \eqref{Pib}. Thus, let 
$Imm_{0}^{m+1,\a}(S^{2}, \bR^{3}) = \Pi^{-1}(\cC^{m,\a} \times [C^{m-1,\a} \setminus \{constants\}])$. 
This is the set of (based) immersions of $S^{2}$ whose images are not spheres of constant mean curvature $c$. 
By the well-known Hopf theorem \cite{H}, these are just the round, constant curvature spheres $S^{2}(r)$. 
One may then consider 
\be \label{Pi2}
\Pi_{0}: Imm_{0}^{m+1,\a}(S^{2}, \bR^{3}) \to \cC^{m,\a} \times [C^{m-1,\a} \setminus \{constants\}],
\ee
$$\Pi_{0}(F) = ([\g], H).$$
For $\Pi_{0}$, conformal reparametrizations of a given immersion act properly on the target data. 

\medskip 

    Alternately, this problem may be remedied in the usual way by choosing a suitable 3-point marking 
for the mapping $F$. Thus, fix three points $p_{i}$ on the round 2-sphere $(\bS^{2}(1), g_{+1}) \subset \bR^{3}$ 
with 
\be \label{n}
dist_{g_{+1}}(p_{i}, p_{j}) = \pi/2,
\ee
for $i \neq j$. Let $Imm_{1}^{m+1,\a}$ be the submanifold of $Imm_{b}^{m+1,\a}$ consisting of 
immersions $F$ such that the pull-back metric $\g = F^{*}(g_{Eucl})$ satisfies \eqref{n}, (with 
$\g$ in place of $g_{+1}$). This is a submanifold of codimension 3, representing a slice to the 
action (by pre-composition) of the conformal group on $Imm_{b}^{m+1,\a}$. Restriction gives an 
induced map
\be \label{slice}
\Pi': Imm_{1}^{m+1,\a}(S^{2}, \bR^{3}) \to \cC^{m,\a} \times C^{m-1,\a}.
\ee
%Funny here that lack of properness is only in $[\g]$ factor and if $H \neq const$, then do have 
%properness of full $\Pi$ map. This all in domain. But in target, doing just the opposite - dividing out 
%by action of conformal group only on $H$ factor, leaving the metric factor fixed. Some kind of 
%duality here. Not fully understood. 

  The map $\Pi'$ is of course still Fredholm. However, while the map $\Pi$ has Fredholm index 0, $\Pi'$ 
now has Fredholm index -3. In order to obtain again a map of index 0, one needs to take a quotient of 
the target space by a 3-dimensional space transverse to the map $\Pi'$. To motivate this, we first examine 
the constraint equations, namely the Gauss-Codazzi and Gauss equations, relating the intrinsic and 
extrinsic geometry of the surface $\S = Im F$. Thus,  for an immersed surface $\S \subset \bR^{3}$, one has: 
\be \label{cg}
\d(A - H\g) = -Ric_{g_{Eucl}}(N, \cdot) = 0,
\ee
\be \label{g}
|A|^{2} - H^{2} + R_{\g} = R_{g_{Eucl}} - 2Ric_{g_{Eucl}}(N,N) = 0.
\ee
Here $A$ is the second fundamental form of the immersion $F$ while $R_{\g}$ and $R_{g_{Eucl}}$ are the scalar 
curvatures of $\g$ and $g_{Eucl}$ respectively. The scalar constraint \eqref{g} will be important later, in Section 3. 

  An immediate consequence of the divergence constraint \eqref{cg} is the following. Suppose $V$ is a 
conformal Killing field on $(S^{2}, \g)$, $\g = F^{*}(g_{Eucl})$. Then pulling back the geometric data 
on $\S$ back to $S^{2}$ via the immersion $F$, one has 
$$\int_{S^{2}}\<V, \d(A - H\g)\>dV_{\g} = \int_{S^{2}}\<\d^{*}V, A - H\g\>V_{\g} = -{\tfrac{1}{2}}
\int_{S^{2}}H div V dV_{\g} = {\tfrac{1}{2}}\int_{S^{2}}V(H)dV_{\g}.$$
Thus \eqref{cg} gives the relation
\be \label{2.6}
\int_{S^{2}}V(H)dV_{\g} = 0,
\ee
for any conformal Killing field $V$ on $(S^{2}, \g)$. In particular, as noted in the Introduction, the relation 
\eqref{2.6} gives an obstruction to the surjectivity of $\Pi$ onto the space of mean curvature functions, 
somewhat analogous to the Kazdan-Warner type obstruction \cite{KW} for prescribed scalar curvature in a 
conformal class. For a given $\g$, \eqref{2.6} is a 3-dimensional restriction on the form of $H$ for 
an immersion $F$. Note however that the constraint \eqref{2.6} is {\it not} a condition on the 
target space data $([\g], H)$, exactly due to presence of the volume form $dV_{\g}$. 

  As in the Introduction, on the target space $C^{m-1,\a}$ of $C^{m-1,\a}$ functions on $S^{2}$, 
define an equivalence relation by 
\be \label{eq2}
[H_{2}] = [H_{1}]  \Leftrightarrow  H_{2} = H_{1} + \ell,
\ee
for some normalized affine function $\ell$ on $S^{2}(1)\subset \bR^{3}$ as following \eqref{eq}. The space of such 
affine functions is exactly $\bR^{3}$, with the origin corresponding to $a = b = 0$ and unit sphere 
corresponding to the unit vector $x \in S^{2}$. This action of $\bR^{3}$ on $C^{m-1,\a}$ is 
free. Let $\cD^{m-1,\a}$ be the space of equivalence classes and let 
\be \label{pi}
\pi: \cC^{m,\a} \times C^{m-1,\a} \to \cC^{m,\a} \times \cD^{m-1,\a},
\ee
be the projection map, (equal to the identity on the first factor). The fibers of $\pi$ are 
$\bR^{3}$, (the space of normalized affine functions). Clearly, the quotient 
$\cC^{m,\a} \times \cD^{m-1,\a}$ is also a smooth Banach manifold. 

  Observe that for any fixed volume form $dV_{\g}$, \eqref{2.6} determines a unique representative 
$H \in [H]$, for the equivalence relation \eqref{eq2}; given any representative $H' \in [H]$, there is a 
unique normalized affine function $\ell$ such that $H' + \ell$ satisfies \eqref{2.6}. 
%?In this context, the 
%submanifold $\cS \subset C$ of functions $H$ satisfying \eqref{2.6} is a smooth slice for the 
%action of the conformal group ${\mathrm{Conf}}(S^{2})$ acting on $C_{1}$ as in \eqref{1.1}. ?
 
  By composing $\Pi'$ in \eqref{slice} with $\pi$ above, we obtain 
\be \label{1.5}
\Pi_{1}: Imm_{1}^{m+1,\a}(S^{2}, \bR^{3}) \to \cC^{m,\a} \times \cD^{m-1,\a},
\ee
$$\Pi_{1}(F) = ([\g], [H]).$$
This is a smooth Fredholm map of Banach manifolds. We claim that the Fredholm index of $\Pi_{1}$ is 
zero, 
\be \label{index}
index(D\Pi_{1}) = 0.
\ee
To prove this, it suffices to consider the linearization $D\Pi_{1}$ at the standard embedding 
$F_{0}: \bS^{2}(1) \subset \bR^{3}$, since the index is invariant under deformations. 
The Fredholm index is 0 if the map $\Pi'$ in \eqref{slice} is transverse to the fibers of $\pi$ in 
\eqref{pi} at the standard embedding $F_{0}$. Thus suppose $X$ is a deformation of the standard
 embedding $F_{0}$. The induced variation of the target data is given by $([\d^{*}X], [H'_{X}])$. To 
 prove transversality, we must show that if $([\d^{*}X], H'_{X}) = (0, \ell)$, i.e.~the variation is 
 tangent to the fiber, then $\ell = 0$, i.e. the variation vanishes in the full target space. 

    To see this, decompose $X = X^{T} + fN$ into tangential and normal components to $\S = Im F$. 
One then has, at $F_{0}$, $2H'_{\d^{*}X} = -\D f - |A|^{2}f + X^{T}(H) = - \D f - 2f$, since $H = |A|^{2} 
= 2$. Suppose then $L(f) = - \D f - 2f = \ell$, for some $\ell$. Now $\ell - a$ is a first eigenfunction 
of the Laplacian on $\bS^{2}(1)$, so that $L(\ell) = L(a) = -2a$. Hence
$$\int_{S^{2}}\ell^{2}dV = \int_{S^{2}}\<L(f), \ell\> = \int_{S^{2}}\<f, L(\ell)\> = -2a\int f.$$
But $\int f = -\frac{1}{2}\int \ell = -2\pi a$, which gives
$$\int_{S^{2}}\ell^{2} = 4\pi a^{2}.$$
However, since $\ell = a + bx$, $\int \ell^{2} = 4\pi a^{2} + b^{2}\int x^{2}$. It follows that $b = 0$ and 
hence $a = 0$, so that $\ell = 0$, proving the claim. 

\medskip 
  
  We will work with both maps $\Pi_{0}$ in \eqref{Pi2} and $\Pi_{1}$ in \eqref{1.5} in the following. 
However, these maps are still not proper, now however for somewhat more subtle reasons than before. 
The basic reason is that immersions may converge to branched immersions maintaining control on 
the data $([\g], H)$. It is worthwhile to first illustrate this clearly on a concrete example. 

  Thus, consider Enneper's surface $E$ in $\bR^{3}$. This is a complete minimal surface immersed in $\bR^{3}$, 
conformally equivalent to the flat plane $\bR^{2}$. The immersion $E: \bR^{2} \to \bR^{3}$ is given explicitly 
by:
$$x_{1} = u(v^{2} - \frac{1}{3}u^{2} + 1),$$
$$x_{2} = -v(u^{2} - \frac{1}{3}v^{2} + 1),$$
$$x_{3} = u^{2} - v^{2}.$$
Equivalently, for $\z = u+iv$, $x_{1} = Re(\z - \frac{1}{3}\z^{3})$, $x_{2} = Im(\z + \frac{1}{3}\z^{3})$, 
$x_{3} = Re(\z^{2})$. We now blow-down $E$ and consider the limit, (the tangent cone at infinity). 
Thus, choose $T$ large and replace $(u, v)$ by $(Tu, Tv)$. Dividing both sides of the equations above 
by $T^{3}$ and relabeling gives
$$x_{1} = u(v^{2} - \frac{1}{3}u^{2} + T^{-2}),$$
$$x_{2} = -v(u^{2} - \frac{1}{3}v^{2} + T^{-2}),$$
$$x_{3} = T^{-1}(u^{2} - v^{2}).$$
Setting $t = T^{-1}$, the curve of blow-downs $E_{t}$ of the Enneper surface is given by
$$x_{1} = u(v^{2} - \frac{1}{3}u^{2} + t^{2}),$$
$$x_{2} = -v(u^{2} - \frac{1}{3}v^{2} + t^{2}),$$
$$x_{3} = t(u^{2} - v^{2}).$$
The immersions $E_{t}$, $t > 0$, are conformal immersions $\bR^{2} \to \bR^{3}$ with 
$H_{E_{t}} = 0$. When $t = 0$, $E_{0}$ is the plane $\bR^{2} = \{x_{3} = 0\}$, with multiplicity 3, 
so a 3-fold branched cover of $\bR^{2}$. Taking the $t$-derivative gives the vector field 
$$X = (2tu, 2tv, u^{2} - v^{2}).$$
The curve of mappings $E_{t}$ is smooth in $t$, for $t \in \bR$. Note that when $t < 0$ one obtains 
a reflection of the ``original" Enneper surfaces $E_{t}$ through the plane $x_{3} = 0$. 

   While the behavior above takes place on the non-compact domain $\bR^{2}$, it can be localized to a 
finite region, then suitably "bent" and extended to a curve of immersions $F_{t}: S^{2} \to \bR^{3}$ with 
$H_{F_{t}} > 0$ everywhere. Note in particular that the immersions $E_{t}$ (or suitably modified 
$F_{t}$) are uniformly bounded, in fact converge, in the $C^{\infty}$ norm and similarly the 
mean curvatures $H_{E_{t}}$ or $H_{F_{t}}$ converge in $C^{\infty}$ as $t \to 0$. In other words, 
the target data $(\g_{t}, H_{t})$ remain bounded (and converge) while the family $E_{t}$ or $F_{t}$ 
does not limit on an immersion. This shows that $\Pi$ is not proper on the space of immersions. 
Of course the Gauss curvature $K$, and the norm of the second fundamental form $A$, blow up as 
$t \to 0$. 

  Similar behavior occurs with the immersions $x_{1} = Re(\z - \frac{1}{2k+1}\z^{2k+1})$, 
$x_{2} = Im(\z + \frac{1}{2k+1}\z^{2k+1})$, $x_{3} = Re(\frac{2}{k+1}\z^{k+1})$, coming from 
the Weierstrass representation of (certain) minimally immersed planes in $\bR^{3}$. Here one has a 
branch point of order $2k+1$ at the origin in the blow-down limit. To obtain branch points of even 
order, in the Weierstrass representation, take $f = z$, $g = z^{k}$. Then $x_{1} = Re(\frac{1}{2}\z^{2} - 
\frac{1}{2k+2}\z^{2k+2})$, $x_{2} = Im(\frac{1}{2}\z^{2} + \frac{1}{2k+2}\z^{2k+2})$, 
$x_{3} = Re(\frac{2}{k+2}\z^{k+2})$. This gives a branch point of even order $2(k+1)$ in the 
blow-down limit. Note however that that this surface is itself a branched minimal immersion; 
setting $w = \z^{2}$ shows that it is a 2-fold branched cover of an immersion. 

   In sum, the construction above can be carried out on any complete minimally immersed plane 
in $\bR^{3}$ of finite total curvature. 

\medskip 

  We now begin a more precise discussion of the space of branched immersions, cf.~\cite{BT}, \cite{GOR}, 
\cite{ET} for further background. Let $S^{2}$ be given its standard two charts, via stereographic projection, 
with local complex coordinate $z = u+iv$ in each chart. This data is fixed throughout the following. 

   Recall that by the uniformization theorem, any immersion $F \in Imm^{m+1,\a}(S^{2}, \bR^{3})$ 
can be reparametrized to a conformal immersion, i.e.~there is a diffeomorphism $\f: S^{2} \to 
S^{2}$ such that $F \circ \f$ is conformal. 

\medskip

{\it Definition}. A map $F \in C^{m+1,\a}(S^{2}, \bR^{3})$ is a conformal branched immersion if it is 
a conformal immersion away from finitely many singular points $\{q_{j}\}$ and in a neighborhood $D$ of 
each singular point $q \in \{q_{j}\}$, the complex gradient $F_{z} = \frac{dF}{dz} =
\frac{1}{2}(F_{u} - iF_{v}): D \to \bC^{3}$ satisfies
\be \label{br}
F_{z} = z^{k}G, 
\ee
in the normalization $z(q) = 0$, where $G: D \to \bC^{3}$ is $C^{m,\a}$ smooth and satisfies
\be \label{con}
G \cdot G = 0, \ \ G(q) \neq 0.
\ee
Also $k \geq 1$; the case $k = 0$ corresponds to a neighborhood where $F$ is an immersion. 

 A map $F \in C^{m+1,\a}(S^{2}, \bR^{3})$ is a branched immersion if it is a reparametrization 
of a conformal branched immersion, so that $F = \w F \circ \f$, where $\w F$ is a conformal 
branched immersion and $\f$ is a diffeomorphism $S^{2} \to S^{2}$. 

\medskip 

  The first condition in \eqref{con} is equivalent to the relation $F_{z}\cdot F_{z} = 0$, i.e.
$$F_{u}\cdot F_{u} - F_{v}\cdot F_{v} - 2iF_{u}\cdot F_{v} = 0,$$
so that $F$ is conformal. The second condition in \eqref{con} then implies that the real and 
imaginary parts of $G$ are linearly independent, and so span the tangent space $T_{F(q)}\S$ 
to the surface $\S = Im F$. 
  
  The mapping $F$ may be recovered from \eqref{br} by solving the inhomogeneous 
Cauchy-Riemann equation, i.e.
$$F(z, \bar z) = \frac{1}{2\pi i}\int_{D}\bar{\xi}^{k}\bar G(\xi, \bar \xi) \frac{dz \wedge d\bar z}{\xi - z} + 
a(\bar z),$$
where $a$ is anti-holomorphic. Since $F$ is real-valued, $a(\bar z)$ is uniquely determined, 
(up to a constant). Observe that elliptic regularity implies that if $G \in C^{m.\a}$ then $F \in C^{m+1,\a}$. 

\medskip 

  Let $BImm^{m+1,\a}(S^{2}, \bR^{3})$ be the set of all branched immersions. This is viewed as 
a subspace of the Banach manifold $Map^{m+1,\a}(S^{2}, \bR^{3})$ and so is given the induced 
topology. 

\begin{lemma} 
The space $BImm^{m+1,\a}(S^{2}, \bR^{3})$ is a smooth Banach manifold. 
\end{lemma}

{\bf Proof:} This is clear for a neighborhood of an immersion $F \in Imm^{m+1,\a}(S^{2}, \bR^{3})$.  
Suppose first $F$ is a conformal branched immersion with singular points $q_{j}$ of the form \eqref{br}. 
Near any such singular point $q$, any conformal branched immersion $\w F$ near $F$ has the form
\be \label{f'}
\w F_{z} = P^{j}(z) \w G^{j} + \nu N, \ \ j = 1, 2, 3,
\ee
where $P^{j}(z)$ is a polynomial of degree $k$ near $z^{k}$, $N$ is the normal vector field to the 
immersion $F$, $\w G$ is near $G$ and $\nu$ is $C^{m,\a}$ smooth with $\nu(q) \sim 0$. 
Of course one requires
$$\w F_{z} \cdot \w F_{z} = 0.$$
For example, if the polynomials $P^{j} = P$ are all equal, with $P = \prod_{i=1}^{\ell}(z - z_{i})^{k_{i}}$, 
then $\w F$ is a nearby conformal branched immersion, with branch points of order $k_{i}$ at $z_{i} 
\sim 0$. Of course $\sum_{1}^{\ell} k_{i} = k$. If on the other hand the three polynomials $P^{j}$ have no 
zeros in common, then $\w F$ is a regular conformal immersion near $F$. (This is the case for instance 
for the curve of Enneper immersions $E_{t}$ above with $t \neq 0$). If $F$ has a branch point of 
order $k$, any {\it sufficiently} near branched immersion has branch points of order $\leq k$; the 
order cannot jump up in sufficiently small neighborhoods. One sees here how curves of 
immersions pass smoothly to branched immersions by allowing distinct roots of the polynomials 
$P^{j}$ to merge together. 

  The tangent space to the space of conformal branched immersions at a conformal branched 
immersion $F$ consists of $C^{m+1,\a}$ vector fields $X$ along $F$. Away from branch points, 
the field $X = dF_{t}/dt$ is $C^{m+1,\a}$ and satisfies the linearization of \eqref{con}, i.e. 
\be \label{conlin}
X_{z}\cdot F_{z} = 0,
\ee
but is otherwise arbitrary. Near a branch point $q$, $X$ has the form 
\be \label{X}
X_{z} = (P^{j})'G^{j} + z^{k}G' + \nu 'N,
\ee
with $G' \in Map^{m,\a}(S^{2}, \bC^{3})$. Since $P_{t}^{j} = \prod_{i=1}^{\ell}(z - z_{i}^{j}(t))^{k_{i}}$, 
one has $(P^{j})' = -z^{k-1}\sum_{i}k_{i}(z_{i}^{j})'$. The linearized equation \eqref{conlin} thus 
becomes  
$$z^{k}[-z^{k-1}\sum_{i}k_{i}(z_{i}^{j})' G^{j}\cdot G + z^{k}G' \cdot G] = 0,$$
near $q$, since $G \cdot N = 0$. Thus, $X$ is tangent to the space of branched conformal 
immersions if and only if  
\be \label{bal}
\sum_{i}k_{i}(z_{i}^{j})' = 0, \ \ {\rm and} \ \ G' \cdot G = 0,
\ee
so that
\be \label{X'}
X_{z} = z^{k}G' + \nu 'N.
\ee
The first condition in \eqref{bal} is a ``balancing condition" on the tangent vectors to the curve of 
roots of $P_{t}^{j}$ at $z = 0$. There is no condition on $\nu'$ besides smoothness. Conversely, 
any $X$ satisfying \eqref{bal} is tangent to a curve $F_{t}$ of conformal branched immersions 
with $F_{0} = F$. This follows from the inverse function theorem, since the map $G' \to G'\cdot G$ 
is surjective onto $C^{m,\a}(D, \bC)$.

   Finally, the space $BImm^{m+1,\a}(S^{2}, \bR^{3})$ consists of arbitrary $C^{m+1,\a}$ smooth 
reparametrizations of conformal branched immersions, i.e.~maps of the form $F\circ \f$, where 
$F$ is a conformal branched immersion and $\f$ is a $C^{m+1,\a}$ diffeomorphism of $S^{2}$. 
Since both factors $F$ and $\f$ have Banach manifold structures, so does the full space of compositions. 
This completes the proof. 
 
{\endproof}

   Note that the space $BImm^{m+1,\a}(S^{2}, \bR^{3})$ is stratified by submanifolds according 
to branching orders. The strata corresponding to branched immersions with total branching order 
$k$ has codimension $6k = dim \, \bC^{3k}$. 

\medskip 

  The mean curvature $H$ of a branched immersion is not well-defined in general at the 
branch points. In fact, for a conformal branched immersion $F$, one has the formula, 
$\frac{1}{4}\D_{0}F = H(F_{u} \times F_{v})$ where $\D_{0} = F_{uu} + F_{vv}$ is the flat 
Laplacian, (cf.~\cite{ET} for instance). In terms of the complex gradient, this is equivalent to 
\be \label{cH}
F_{\bar z z} = iH(\bar{F_{z}}\times F_{z}). 
\ee
One has 
$$\bar{F_{z}}\times F_{z} = |z|^{2k} \bar G \times G,$$
and so for the mean curvature to be well-defined in $C^{m-1,\a}$ one needs
\be \label{cH'}
\frac{F_{\bar z z}}{|z|^{2k}} \in C^{m-1,\a},
\ee
in a neighborhood of a branch point of order $k$. (Of course this holds automatically away 
from branch points). Conversely, by elliptic regularity associated to \eqref{cH}, if $H \in C^{m-1,\a}$, 
then $F \in C^{m+1,\a}$. Since $F_{z} = z^{k}G$, so $F_{\bar z z} = z^{k}G_{\bar z}$ it is more 
convenient to express the equation above as
$$G_{\bar z} = \bar z^{k}\f N,$$
for some $\f \in C^{m-1,\a}(S^{2}, \bC^{3})$, where $N$ is a unit normal for the surface 
$\S = Im F$. Note that $\bar F_{z}\times F_{z}$ and $\bar G \times G$ are multiples of $N$. 

\medskip 

{\it Definition}. A $C^{m+1,\a}$ conformal branched immersion $F$ is $H$-regular if near any branch 
point $q$ of order $k$ one has
\be \label{G}
G_{\bar z} = \bar z^{k}\f N, 
\ee
for some $\f \in C^{m-1,\a}(D, \bC)$. A $C^{m+1,\a}$ branched immersion $F$ is 
$H$-regular if it is a reparametrization (by a $C^{m+1,\a}$ diffeomorphism) of a 
conformal branched $H$-regular immersion. 

\medskip 

   The mean curvature $H_{F}$ of an $H$-regular branched immersion is well-defined, and in 
$C^{m-1,\a}(S^{2})$. The space of $H$-regular branched immersions is denoted by 
$\cH Imm^{m+1,\a}(S^{2}, \bR^{3})$ and is a closed submanifold of $BImm^{m+1,\a}(S^{2}, \bR^{3})$. 
In particular $\cH Imm^{m+1,\a}(S^{2}, \bR^{3})$ is itself a smooth Banach manifold. 

  The linearization of the mean curvature $H$ among conformal branched immersions is given 
by 
\be \label{Hlin}
G'_{\bar z} = \bar z^{k}(iH'_{X}\bar G \times G + iH(\bar G'\times G + \bar G \times G')),
\ee
so that the linearization of \eqref{G} holds with $\f' \in C^{m-1,\a}$.

\medskip 
  
   Next we also need to enlarge the target manifold $\cC^{m,\a}$. First, the metric $\g$ 
on $S^{2}$ is written in terms of the coordinates $(z, \bar z)$ in place of 
the real and imaginary parts of $z$. In the following, its convenient to let $x = u + iv = 
u_{1} + iu_{2}$. Thus 
\be \label{cm}
\g = \g_{zz}dz^{2} + 2\g_{\bar z z}d\bar z dz + \g_{\bar z \bar z}d\bar z^{2}.
\ee
In terms of real coordinates $\g = \g_{ij}du_{i}du_{j}$, one has the relations 
$$\g_{zz} = {\tfrac{1}{4}}(\g_{11}-\g_{22} - 2i\g_{12}), \ \ \g_{\bar z z} = 
{\tfrac{1}{2}}(\g_{11} + \g_{22}), \ \ \g_{\bar z \bar z} = \bar \g_{zz}.$$

   Now let $Met_{s}^{m,\a}(S^{2})$ be the space of $C^{m,\a}$ symmetric bilinear forms $\g$ on 
$S^{2}$ which are positive definite outside a finite number of singular points $q_{j}$ and which 
near each singular point $q \in \{q_{j}\}$, have the form \eqref{cm} with 
$$\g_{zz} = z^{2k}\f_{1}, \ \ \g_{\bar z z} = |z|^{2k}\f_{2},$$
where $\f_{1}\in C^{m,\a}(D, \bC)$ and $\f_{2} \in C^{m,\a}(D, \bR^{+})$ near $q$. 

   Clearly $Met_{s}^{m,\a}(S^{2})$ is a smooth Banach manifold. The tangent space consists of 
$C^{m,\a}$ symmetric bilinear forms $h$ of the form \eqref{cm} (no longer necessarily positive 
definite) with 
\be \label{mv}
h_{zz} = z^{2k}\f_{1}', \ \ h_{\bar z z} = |z|^{2k}\f_{2}'.
\ee
   
   Define an equivalence relation on 
$Met_{s}^{m,\a}(S^{2})$ by setting $\g_{2} \sim \g_{1}$ if $\g_{2} = \mu^{2}\g_{1}$, for some 
positive function $\mu \in C^{m,\a}(S^{2})$. Abusing notation, let $\cC^{m,\a}$ denote the 
quotient space. This is the space of pointwise conformal equivalence classes of singular, branched 
metrics on $S^{2}$.

   Now given an $H$-regular branched immersion $F \in \cH Imm^{m+1,\a}(S^{2}, \bR^{3})$, 
the induced or pullback metric $F^{*}(g_{Eucl})$ is a singular branched metric on $S^{2}$, so gives 
an element in $Met_{s}^{m,\a}(S^{2})$. Similarly, the mean curvature $H_{F}$ is a $C^{m-1,\a}$ smooth 
function on $S^{2}$. Thus, one has a map 
\be \label{pi0}
\Pi_{0}: \cH Imm^{m+1,\a}(S^{2}, \bR^{3}) \to \cC^{m,\a} \times [C^{m-1,\a}(S^{2}) 
\setminus \{constants\}],
\ee
$$\Pi_{0}(F) = ([F^{*}(g_{Eucl})], H_{F}),$$ 
extending the map $\Pi_{0}$ in \eqref{Pi2}. From its construction, $\Pi_{0}$ is a smooth map 
of Banach manifolds. Similarly, the map
\be \label{pi1}
\Pi_{1}: \cH Imm^{m+1,\a}(S^{2}, \bR^{3}) \to \cC^{m,\a} \times \cD^{m-1,\a},
\ee
$$\Pi_{1}(F) = ([F^{*}(g_{Eucl})], [H_{F}]),$$ 
extending the map $\Pi_{1}$ in \eqref{1.5} is a smooth map of Banach manifolds.

\begin{proposition}
The map $\Pi_{0}$ in \eqref{pi0} is Fredholm, of Fredholm index 0.  
\end{proposition}

{\bf Proof:} This is proved in \cite{An1}, \cite{An2} when $\Pi_{0}$ is restricted to regular 
immersions, as in \eqref{Pi2}. It is useful present the short proof here. 

   The linearization $D\Pi_{0}$ acts on vector fields $X$ along the immersion $F$. Write 
$X = X^{T} + fN$, where $X^{T}$ is tangent and $N$ is normal to $\S = Im(F)$. Then 
\be \label{3.4}
\d^{*}X = \d^{*}(X^{T}) + fA + df\cdot N,
\ee
so that $(\d^{*}X)^{T} = \d^{*}(X^{T}) + fA$. The second term here is lower order in $X$ 
and so does not contribute to the principal symbol. The principal symbol $\sigma$ of 
the first component of $D\Pi_{0}$ is thus
\be \label{3.5}
\sigma([2(\d^{*}X)^{T}]_{0}) = \sigma([2(\d^{*}(X^{T})]_{0}) = 2(\xi_{i}X_{j} - \frac{\xi_{i}X_{i}}{2}\delta_{ij}),
\ee
where $i,j$ are indices along $\S$, i.e.~tangent to $S^{2}$. Setting this to 0 gives
$$\xi_{1}X_{2} = \xi_{2}X_{1} = 0 \ \ {\rm and} \ \ \xi_{1}X_{1} = \xi_{2}X_{2}.$$
Since $(\xi_{1}, \xi_{2}) \neq (0,0)$, it is elementary to see that the only solution 
of these equations is $X_{1} = X_{2} = 0$. Next, for the mean curvature, one has 
$2H_{\delta^{*}X}' = -\D f - |A|^{2}f + X^{T}(H)$. Here $\D$ is the Laplacian with respect 
to the induced metric $\g = F^{*}(g_{Eucl})$. Hence $\D = \l^{-2}\D_{0}$, for $\g = \l^{2}g_{0}$.  
The leading order symbol acting on $f$ is thus $\l^{-2} |\xi|^{2}f$. For a regular immersion, 
$\l > 0$ and the vanishing of this term implies $f = 0$. Thus, the symbol of 
$D\Pi_{0}$ is elliptic, so that by the regularity theory for elliptic systems, cf.~\cite{Mo} 
for instance, $D\Pi_{0}$ is Fredholm. Note that the operator $D\Pi_{0}$ has terms 
containing both first and second derivatives of $X$. 

   The proof at or near branch points of $F$ is very similar, and we carry out the details 
below. Without loss of generality, we may assume that $F$ itself is conformal. 

For $F_{z} = \frac{1}{2}(F_{u} - iF_{v})$ and $F_{\bar z} = \bar{F_{z}} = \frac{1}{2}
(F_{u} + iF_{v})$ one has
$$4F_{z} \cdot F_{z} = F_{u}\cdot F_{u} - F_{v}\cdot F_{v} - 2iF_{u}\cdot F_{v},$$
so that $F$ conformal if and only if 
$$F_{z}\cdot F_{z} = 0.$$
Similarly
$$4F_{z}\cdot F_{\bar z} = F_{u}\cdot F_{u} + F_{v}\cdot F_{v} \equiv \l^{2},$$
with $\l^{2} \in C^{m,\a}(S^{2}, \bR)$. Thus, the non-conformal variation of the 
metric is given by 
$$\frac{d}{dt}F_{z}^{t}\cdot F_{z}^{t} = 2X_{z}\cdot F_{z},$$
where $X = dF_{t}/dt$. From \eqref{X'} one then has 
\be \label{F1}
X_{z}\cdot F_{z} = z^{2k}G'\cdot G.
\ee
This is exactly of the form \eqref{mv} with 
$$\f_{1}' = G' \cdot G.$$
One may choose a basis of $\bR^{3}$ so that $G^{3} = 0$ so that the tangent space to 
$\S = Im F$ at $F(q) = 0$ is the $(x_{1}, x_{2})$ plane $\bR^{2} \subset \bR^{3}$. Then 
\eqref{F1} is clearly Fredholm in $(G^{j})'$ for $j = 1,2$ and so Fredholm in the first two 
components of $X_{z}$, (just as above following \eqref{3.5}). 

   Next, as in \eqref{Hlin}, the linearization $H'$ of the mean curvature is given by 
$$G'_{\bar z} = \bar z^{k}(iH'_{X}\bar G \times G + iH(\bar G'\times G + \bar G \times G')).$$
As before, given $(G^{j})'$, $j = 1, 2$, this gives an elliptic equation for the remaining third 
component $(G^{3})'$ of $G'$, or equivalently an elliptic equation for the third component of 
$X_{z}$. It follows that $D\Pi_{0}$ is Fredholm. By the deformation invariance of the index, 
$D\Pi_{0}$ is of Fredholm index 0. 

\endproof

\begin{proposition}
The map $\Pi_{1}$ in \eqref{pi1} is Fredholm, of Fredholm index 0.  
\end{proposition}

{\bf Proof:} This follows directly from Proposition 2.2 and \eqref{index}. 

{\endproof}

\section{Properness}
  
  In this section, we prove that the maps $\Pi_{0}$ and $\Pi_{1}$ in \eqref{pi0} and \eqref{pi1} 
are proper, when restricted to the domains where $H_{F} > 0$. Given Propositions 2.2 and 2.3, 
this is the key to the proof of Theorem 1.1. 

   Proving that $\Pi_{0}$ is proper requires showing that an immersion (or branched immersion) 
$F: S^{2} \to \bR^{3}$ is controlled in $C^{m+1,\a}$ by the target data $([F^{*}(g_{Eucl}], H_{F})$. 
The starting point to obtain such control is the Gauss constraint \eqref{g}. Thus, integrating 
\eqref{g} along the immersion $F$ and using the Gauss-Bonnet theorem gives 
\be \label{2.15}
\int_{S^{2}}|A|^{2} = \int_{S^{2}}H^{2} - 4\pi \chi(S^{2}) = \int_{S^{2}}H^{2} - 8\pi.
\ee
Since $H$ is pointwise controlled by the target data, \eqref{2.15} gives control on the $L^{2}$ 
norm of $A$, provided one can obtain a bound on the area of $\S = Im F$. Note this does not 
follow directly from control of the conformal class $[\g]$. 

   Let $F_{t}$ be a curve in $\cH Imm^{m+1,\a}(S^{2}, \bR^{3})$ with $F_{0}$ a fixed embedding - say 
near the standard round sphere. Although the mean curvature $H$ controls the variation of the area of 
$\S_{t} = Im F_{t}$ in the normal direction, pointwise control of $H$ alone does not give rise to an area 
bound. For example, consider a long cylinder $C = [-L, L]\times S^{1}(1)$ embedded in $\bR^{3}$, of 
mean curvature $H = 1$. One may attach two spherical caps to the boundary $\partial C$ to obtain an 
embedding $F_{L}: S^{2} \subset \bR^{3}$ with $1 \leq H \leq 2$ pointwise. As $L \to \infty$, 
$$area (\S_{L}) \to \infty,$$
with uniform control on $H$; here $\S_{L} = Im F_{L}$. Performing a similar construction with the 
periodic Delaunay cylinders of constant mean curvature, one may find curves of immersions 
$F_{t}$ with the same behavior with 
$$1 \leq H \leq 1+\e,$$
for any fixed $\e > 0$. It is important to note however that the (pointwise) conformal structure 
of $(S^{2}, \g_{t})$, $\g_{t} = F_{t}^{*}(g_{Eucl})$ (or $\g_{L} = F_{L}^{*}(g_{Eucl})$) degenerates, 
i.e.~diverges to infinity, in the examples above. In particular, the mean curvature viewed as a 
function $H: S^{2} \to \bR$ with a fixed (standard) atlas on $S^{2}$ does not remain uniformly 
bounded in $C^{m-1,\a}$ in these examples. 
  
\medskip 
  
  We first prove the area estimate for $\Pi_{0}$, and then show a similar argument gives the 
result for $\Pi_{1}$. Let $\cH Imm_{+}^{m+1,\a} = \cH Imm_{+}^{m+1,\a}(S^{2}, \bR^{3})$ be 
the space of $H$-regular branched immersions with mean curvature $H_{F} > 0$. 

\begin{proposition}
For $F \in \cH Imm_{+}^{m+1,\a}(S^{2}, \bR^{3})$, suppose the target data $([\g], H) = 
([F^{*}(g_{Eucl})], H_{F})$ satisfy
\be \label{bound}
||([\g], H)|| \leq L,
\ee
where the norm is taken in the normed target space $\cC^{m,\a} \times C_{+}^{m-1,\a}$. Suppose 
also 
\be \label{low}
||H - c||_{C^{m-1,\a}} \geq L^{-1},
\ee
for any constant $c$.

   Then there is a constant $A_{0}$, depending only on $L$, such that 
\be \label{A}
area(F) \leq A_{0}.
\ee
\end{proposition}

{\bf Proof:}  Since any branched immersion can be perturbed slightly to an immersion with a 
correspondingly small effect on the target data $([\g], H)$, it suffices to prove the result for 
$F \in Imm_{+}^{m+1,\a}(S^{2}, \bR^{3})$.  By a well-known result of Smale \cite{Sm1}, the space 
$Imm^{m+1,\a}(S^{2}, \bR^{3})$ is connected. Moreover, it then follows from \cite{LM} that 
$Imm_{+}^{m+1,\a}(S^{2}, \bR^{3})$ is connected. 

  Fix an embedding $F_{0}: S^{2} \to \bR^{3}$ near that of the standard round sphere. 
Let $F_{t}$, $0 \leq t \leq L$ be a curve in $Imm_{+}^{m+1,\a}$ starting at $F_{0}$ and ending 
at $F_{L} = F$. We may assume that $F_{t}$ is parametrized in such a way that the image 
curve $\Pi(F_{t})$ is parametrized by (or proportional to) arclength, so that 
\eqref{bound}-\eqref{low} hold along $F_{t}$. 

  Let $X = dF_{t}/dt$ be the variation vector field of the curve $F_{t}$, and let 
$\k = (\cL_{X}g)^{T} = 2(\d^{*}X)^{T}$ be the variation of the induced metric $\g$.
Write $X = X^{T} + fN$, so that $\frac{1}{2}\k = \d^{*}X^{T} + fA$. By the standard second 
variation formula for area, the variation of the mean curvature is given by 
\be \label{HX}
H'_{X} = -\D f - |A|^{2}f + X^{T}(H).
\ee
Integrating this over $S^{2}$ gives
\be \label{H'}
\int_{S^{2}}H'_{X} = \int_{S^{2}}-|A|^{2}f + X^{T}(H).
\ee

   By the uniformization theorem, any metric $\g$ on $S^{2}$ is of the form $\g = 
\l^{2}\psi^{*}(\g_{+1})$ for some function $\l > 0$ and diffeomorphism $\psi$ of $S^{2}$. By 
\eqref{bound}, the conformal class of $[\g]$ is uniformly bounded and hence the diffeomorphism 
$\psi$ is uniformly controlled in $C^{m+1,\a}$ (by $L$) modulo the conformal group.  As discussed 
in Section 2, since $H$ is bounded away from the constant functions, the control on $H$ gives control 
on the conformal group factor. (The conformal group acts properly on the target data $([\g], H)$ 
when $H$ is non-constant). It follows that one may precompose $F_{t}$ with a bounded curve of 
diffeomorphisms $\psi_{t}$ so that the conformal class is fixed, i.e.~$F_{t}^{*}g_{Eucl} = 
\l_{t}^{2}g_{+1}$. This gives
\be \label{k}
{\tfrac{1}{2}}\k = \d^{*}X^{T} + fA = \f \g,
\ee
for some (undetermined) conformal factor $\f = \f_{t}$. The diffeomorphisms $\psi_{t}$ alter 
$H_{F_{t}}$ and $H'_{X}$ only by a uniformly bounded factor, which is ignored in the following. 

  Now compute:
$$(\int_{S^{2}}H_{F_{t}}dV_{\g_{t}})' = \int_{S^{2}}H'_{X} + \int_{S^{2}}H{\tfrac{1}{2}}tr \k = 
\int_{S^{2}}H'_{X} + 2H\f.$$
Pairing \eqref{k} with $A$ gives
$$\<A, \d^{*}X^{T}\> + f|A|^{2} = \f H,$$
so that 
$$\int_{S^{2}}\f H =  \int_{S^{2}}f|A|^{2} + \<X^{T}, \d A\>.$$
By the divergence constraint \eqref{cg}, $\d A = \d(H\g) = -dH$, and hence
$$\int_{S^{2}}2\f H = 2\int_{S^{2}}f|A|^{2} - X^{T}(H)  = -2\int_{S^{2}}H'_{X},$$
where the last equality follows from \eqref{H'}. In sum, for such conformal variations, one has 
\be \label{Hest}
(\int_{S^{2}}H)' = -\int_{S^{2}}H'_{X},
\ee
so that 
\be \label{Hest2}
|(\int_{S^{2}}H_{F_{t}}dV_{\g_{t}})'| \leq  K area(F_{t}),
\ee
where $K$ is a bound for $|H'_{X}|$, (cf.~the statement following \eqref{k}). Since $H$ is 
uniformly controlled, 
\be \label{Hcon}
0 < H_{0} \leq H \leq H_{0}^{-1},
\ee
integrating over $t$ gives
$$H_{0}area(F_{t}) \leq K \int_{0}^{t}area(F_{t}) + c,$$
which is the same as the differential inequality $H_{0}f' \leq Kf + c$, for $f = \int area(F_{t})$. It follows 
by a simple calculus argument that 
$$area(F_{t}) \leq Ce^{K_{1}t} \leq Ce^{K_{1}L}.$$

{\endproof}  

   Next we prove the analog of Proposition 3.1 for the map $\Pi_{1}$. 

\begin{proposition}
For $F \in \cH Imm_{+}^{m+1,\a}$, suppose the target data $([\g], [H]) = ([F^{*}(g_{Eucl})], [H_{F}])$ 
for $\Pi_{1}$ satisfy
\be \label{con2}
||([\g], [H])|| \leq L,
\ee
where the norm is taken in the normed target space $\cC^{m,\a} \times \cD_{+}^{m-1,\a}$. 
Then there is a constant $A_{0}$, depending only on $L$, such that 
\be \label{A1}
area(F) \leq A_{0}.
\ee
\end{proposition}

{\bf Proof:}  The proof is essentially the same as that of Proposition 3.1. In this case, the control over 
the diffeomorphisms $\psi_{t}$ comes from the normalization \eqref{n} giving a slice to the action 
of the conformal group $\mathrm{Conf}(S^{2})$. One also needs to show that control over $[H_{F_{t}}]$ 
implies control over $H_{F_{t}}$. To see this, write $H_{F_{t}} = H_{t} + \ell_{t}$ where $H_{t}$ is uniformly 
controlled (by $L$) in $C^{m-1,\a}$. Let  $X_{t}$ be the conformal vector field generated by the linear function 
$x_{t}$, where $\ell_{t} = a_{t} + b_{t}x_{t}$. Then by \eqref{ob}, 
\be \label{ellcont}
\int_{S^{2}}(X_{t}(H_{t}) + b_{t}\f)dV_{\g_{t}} = 0,
\ee
where $\f = X_{t}(x_{t}) = \sin r$, $0 \leq r \leq \pi$, is a fixed function on $S^{2}$ up to rotation, 
$0 \leq \f \leq 1$, which vanishes only at the poles. Since $X_{t}(H_{t})$ is uniformly controlled, 
it follows that the family $b_{t}$ is uniformly bounded, so that $\ell_{t}$ is also uniformly bounded, 
which proves the claim. Finally, since $\min(H+\ell) \geq \min H$ for all normalized affine functions, 
\eqref{Hcon} remains valid under the assumption \eqref{con2}. 

{\endproof}  
  
  The area bounds \eqref{A} or \eqref{A1} together with \eqref{2.15} give an apriori bound on the 
scale-invariant quantity $\int |A|^{2}$,
\be \label{2.16}
\int_{S^{2}} |A|^{2} \leq C,
\ee
for $H$-regular branched immersions, with $C$ depending only on the norm of the target data for 
$\Pi_{0}(F)$ or $\Pi_{1}(F)$.

\begin{theorem} 
The maps 
\be \label{pi01}
\Pi_{0}: \cH Imm_{+}^{m+1,\a}(S^{2}, \bR^{3}) \to \cC^{m,\a} \times [C^{m-1,\a}(S^{2}) 
\setminus \{constants\}],
\ee
$$\Pi_{0}(F) = ([F^{*}(g_{Eucl})], H_{F}),$$ 
and
\be \label{pi11}
\Pi_{1}: \cH Imm_{+}^{m+1,\a}(S^{2}, \bR^{3}) \to \cC^{m,\a} \times \cD^{m-1,\a},
\ee
$$\Pi_{1}(F) = ([F^{*}(g_{Eucl})], [H_{F}]),$$ 
i.e.~the maps in \eqref{pi0}-\eqref{pi1} restricted to $\cH Imm_{+}^{m+1,\a}(S^{2}, \bR^{3})$, 
are smooth proper Fredholm maps.
\end{theorem}

{\bf Proof:} We first work with $\Pi_{0}$. Suppose $F_{i}$ is a sequence in $\cH Imm_{+}^{m+1,\a}$ 
such that the target data $\Pi_{0}(F_{i}) = ([\g_{i}], H_{i})$ converge in $\cC^{m,\a}\times C^{m-1,\a}$ 
to a limit $([\g], H) \in \cH Imm_{+}^{m+1,\a}$. One then needs to show that a subsequence of $\{F_{i}\}$ 
converges in $C^{m+1,\a}$ to a limit map $F \in \cH Imm_{+}^{m+1,\a}$. As in the proof of Proposition 3.1, 
for simplicity we may assume that the maps $F_{i}$ are regular immersions. 

   It is standard and well-known that if there is a uniform bound for the second fundamental 
form $A = A_{F_{i}}$ of $\{F_{i}\}$, then a subsequence converges in $C^{1,\a}$ to a limit 
$C^{1,\a}$ immersion $F$. (The bound on $|A|$, together with control of the diffeomorphisms 
reparametrizing $F_{i}$ as in the proof of Proposition 3.1, imply a bound on the second derivatives of 
the immersion $F$; the result then follows from the Arzela-Ascoli theorem). Since the 
data $([F^{*}g_{Eucl}], H)$ are elliptic for the map $F$, elliptic regularity then shows that the 
convergence of the target data implies convergence $F_{i} \to F$ in $C^{m+1,\a}$. 

   Thus the issue is to understand the structure of $\{F_{i}\}$ when $|A_{i}|$ blows up as 
$i \to \infty$. This is done by a blow-up argument. Choose a point $x = x_{i}$ on $S^{2}$ 
where $|A_{F_{i}}|$ is (locally) maximal, so that if $|A_{F_{i}}|(x_{i}) = \l_{i}$, then 
$\l_{i} \to \infty$ as $i \to \infty$ with 
$$|A_{F_{i}}|(y) \leq \l_{i},$$
with $y$ near $x_{i}$. Now rescale the immersion at $x_{i}$ by multiplying by 
$\l_{i}$, i.e.~consider $F^{\l_{i}}_{i} = \l_{i} F_{i}$. Here we assume (without loss of generality) 
that $F_{i}(x_{i}) = 0 \in \bR^{3}$. We also rescale or blow-up the local coordinates for 
$S^{2}$ near $x_{i}$ by $\l_{i}$, exactly as in the discussion of the Enneper surfaces in 
Section 2. Note that norms of derivatives of $F_{i}$ are invariant under such simultaneous 
rescalings of domain and range. One then has 
\be \label{n1}
|A_{F^{\l_{i}}_{i}}|(y_{i}) \leq 1, \\ {\rm with} \ \ |A_{F^{\l_{i}}_{i}}|(x_{i}) = 1.
\ee
This holds for all $y_{i}$ such that $F_{i}^{\l_{i}}(y_{i})$ is of uniformly bounded distance to 
$F_{i}^{\l_{i}}(x_{i}) = 0$. It follows from the constraint equation \eqref{g} that the Gauss curvature 
of the blow-up surfaces $F_{i}^{\l_{i}}$ remains uniformly bounded. 

   The standard compactness result used above thus implies that 
\be \label{locconv}
F_{i}^{\l_{i}} \to F_{\infty} \ \ {\rm in} \ \ C_{loc}^{m+1,\a},
\ee
and $F_{\infty}$ is an immersion $F_{\infty}: \bR^{2} \to \bR^{3}$ with $x_{i} \to x_{\infty}$ and 
$F_{\infty}(x_{\infty}) = 0$. 
Moreover, since 
$$H_{i}^{\l} = \l_{i}^{-1}H_{F_{i}} \to 0 \ \ {\rm as} \ \ i \to \infty,$$
$F_{\infty}$ is a complete minimal immersion $\bR^{2} \to \bR^{3}$. Also, since the conformal 
structure of $F_{i}$ is uniformly controlled, $[F_{\infty}^{*}(g_{Eucl}] = [du^{2} + dv^{2}]$ is the 
standard conformal structure on $\bR^{2}$. 

  The smooth convergence in \eqref{locconv} implies that the limit immersion $F_{\infty}$ is not 
totally geodesic, since
\be \label{n2}
|A_{F_{\infty}}|(x_{\infty}) = 1.
\ee
Moreover, since the $L^{2}$ norm of $A$ is scale invariant, the bound \eqref{2.16} holds uniformly 
for the family $F_{i}^{\l_{i}}$ and hence the limit minimal immersion $F_{\infty}$ has finite total curvature, 
$$\int_{\bR^{2}}|A|_{F_{\infty}}^{2} < \infty.$$

  By a well-known result of Osserman, cf.~\cite{O} for example, the total scalar curvature 
$R = -|A|^{2}$ of a minimal surface $\S$ immersed in $\bR^{3}$ is quantized, i.e.
$$\int_{\Sigma}|A|^{2} = 4k\pi,$$
with $k = 0$ exactly when $\Sigma$ is a flat totally geodesic plane $\bR^{2} \subset \bR^{3}$ 
while $k = 1$ exactly for Enneper's surface, (up to scaling). 

   It follows directly from \eqref{n2} that $k \geq 1$, on any sequence $x_{i}$ as above in \eqref{n1}. 
In addition, the sequence $F_{i}$ satisfies the scale-invariant uniform bound \eqref{2.16}. 
Hence there at most 
$$N \leq \frac{C}{4\pi},$$
points $q_{j}$ (limit points of sequences $\{x_{i}\}$) where $|A|$ can blow up. 

   Away from those points, one has smooth convergence in $C^{m+1,\a}$ to a limit immersion $F$ - 
as discussed above in the second paragraph of the proof. The domain of $F$ here is a finitely punctured 
two sphere $S^{2} \setminus \cup \{q_{j}\}$, $1 \leq j \leq N$. The singular points correspond to formation 
of  branch points at $F$. In more detail, near such points $x_{i}$, $F_{i}$ is $C^{m+1,\a}$ close to the 
the $\l_{i}^{-1}$-blow-down $E_{\l_{i}^{-1}}$ of a complete minimally and conformally immersed plane 
$E \simeq \bR^{2} \to \bR^{3}$. Recall as above that the $C^{m+1,\a}$ norm of $F_{i}$ is invariant under 
the rescalings above. The limit $E_{0}$ of the minimal immersions $E_{\l_{i}^{-1}}$ is the map $h(z) = 
(z^{k}, 0): \bC \to \bC\times \bR = \bR^{3}$, i.e.~a standard branched immersion with branch point 
of order $k-1 \geq 2$. From this, it is readily verified that $F \in \cH Imm^{m+1,\a}$ and 
$F_{i} \to F$ in $\cH Imm^{m+1,\a}$. This proves that $\Pi$ is proper. 

   To prove that $\Pi_{1}$ is proper, it suffices from the above to prove that if a sequence of immersions 
$F_{i}$ satisfies $[H_{F_{i}}| \to [H]$ then $H_{F_{i}}$ converges (in a subsequence) in $C^{m-1,\a}$. 
The proof of this is the same as the proof of Proposition 3.2, i.e.~as in \eqref{ellcont}, with $t$ 
replaced by $i$. 

{\endproof} 

  A proper Fredholm map $\Pi: \cB_{1} \to \cB_{2}$ of index 0 between Banach manifolds has 
a well-defined degree (mod 2), (the Smale degree) given by
\be \label{deg}
deg \, \Pi = \# \Pi^{-1}(y), \ \ ({\rm mod} \ 2),
\ee
for any regular value $y \in \cB_{2}$. The regular values of $\Pi$ are open and dense in $\cB_{2}$ 
and the properness of the map $\Pi$ ensures that the cardinality in \eqref{deg} is finite. In many 
situations, the $\bZ_{2}$-valued degree can be enhanced to a $\bZ$-valued degree; one needs 
suitable orientations for the triple $\Pi: \cB_{1} \to \cB_{2}$, cf.~\cite{ETr} for example. However, we 
will not explore this further here. 

\begin{proposition}
For the maps $\Pi_{0}$ and $\Pi_{1}$ in \eqref{pi01} and \eqref{pi11}, one has 
\be \label{d0}
deg \, \Pi_{0} = 0,
\ee
but
\be \label{d1}
deg \, \Pi_{1} = 1.
\ee
\end{proposition}  
  
{\bf Proof:}  It is well-known (and easy to see) that any proper Fredholm map of index 0 between Banach 
manifolds which is not surjective has degree 0. Clearly $\Pi_{0}$ cannot be surjective, exactly 
due to the obstruction \eqref{ob}. This gives \eqref{d0}. 

   To determine $deg \, \Pi_{1}$, consider the value $([\g_{+1}], [2])$, so that $H = 2 + \ell$ for some 
affine function $\ell$ with $|\ell| \leq 2$. This corresponds to data for the standard embedding $\bS^{2}(1) 
\subset \bR^{3}$. Clearly the obstruction \eqref{ob} is satisfied for the functions $H = 2 + \ell$ only 
for $\ell = 0$. By a classical theorem of Hopf, any immersed sphere of constant mean curvature 2 is 
a reparametrization of the standard embedding. Fixing the conformal class to be $[\g_{+1}]$ implies that 
the parametrization is conformal, and the 3-point condition \eqref{n} implies that the mapping is the 
standard embedding.  It follows that the standard embedding uniquely realizes this value for $\Pi_{1}$. 

   We claim that $Ker D\Pi_{1} = 0$ at the standard embedding $F_{0}$. Any $k \in Ker D\Pi_{1}$ satisfies
\be \label{ker}
k^{T}_{0} = 0, \ \ H'_{k} = \ell,
\ee
for some $\ell$. As in the proof of Proposition 3.1, $k = 2\d^{*}Z$, for some vector field $Z$ along $F_{0}$. Write 
$Z = Z^{T} + fN$, where $Z^{T}$ is tangent to $\bS^{2}(1)$ and $N$ is the unit outward normal. Then  
$\frac{1}{2}k = \d^{*}Z = \d^{*}Z^{T} + fA$. Since $A = \g$, the first equation in \eqref{ker} gives
$$\d^{*}Z^{T} = \f \g,$$
with $\f = \frac{1}{2}div Z$. Thus $Z^{T}$ is a conformal Killing field on $\bS^{2}(1)$. The $3$-point 
normalization \eqref{n} then forces $Z^{T} = 0$. The vector field $fN$ is a ``Jacobi field" along 
$\bS^{2}(1)$, i.e.~satisfies 
$$\D f + |A|^{2}f = \D f + 2f = \ell.$$
The same argument as following \eqref{index} shows that necessarily $\ell = 0$. Hence $f$ is a  
first eigenfunction of the Laplacian on $\bS^{2}(1)$, corresponding to the normal components of Killing fields 
$T$ (translations) on $\bS^{2}(1)$. However, such infinitesimal translations violate the normalization \eqref{Pib} 
that the immersions $F$ are based immersions. Thus $k =0$ and so $Ker D\Pi_{1} = 0$. 
It follows that $deg \, \Pi_{1} = 1$, which proves the result. 

{\endproof}

   The fact that $deg \, \Pi_{1} = 1$ implies that $\Pi_{1}$ is surjective, which is just the statement of 
Theorem 1.1. Thus Theorem 1.1 is proved. 

\begin{remark}
{\rm The contrast of the two degrees in \eqref{d0} and \eqref{d1} is rather unusual. Theorem 1.1 
gives the existence of a branched immersion $F: S^{2} \to \bR^{3}$ realizing any 
prescribed $([\g], [H])$, so that $H_{F} = H + \ell$, where $H > 0$ is arbitrarily prescribed. 
On the other hand, since $deg \, \Pi_{0} = 0$, given one such immersion $F$ with data $([\g], H_{F})$, 
there must exist generically at least one more distinct immersion $F'$, giving at least two branched immersions 
realizing $([\g], H_{F})$. Here generic means $([\g], H_{F})$ is a regular value of $\Pi_{0}$. 

   As a concrete example, it follows that for $\e$ sufficiently small, the data $([\g], H)$ with 
$0 < |H - c| < \e$ near the standard round spherical data are realized by at least two distinct 
immersions $F$, $F'$. 

   Such non-congruent pairs of immersions $F$ with equal values of $([\g], H)$ may be considered 
as ``conformal Bonnet pairs". Recall that a Bonnet pair is a pair of immersions which are isometric 
and with identical mean curvatures. It is well-known that there are no Bonnet pairs of immersions 
$S^{2} \to \bR^{3}$, cf.~\cite{LT}, \cite{Sa}.  

}
\end{remark}

\section{Generalizations}

  In this section, we discuss several generalizations of Theorem 1.1. 
  
   First, many of the results above apply to surfaces $\S$ of genus $g > 0$. In this case, there is no 
obstruction to the form of the mean curvature $H$ as in \eqref{ob}. This is immediate when $g > 1$ 
since such surfaces have no conformal vector fields. In the case of $g = 1$, there are conformal 
vector fields on a torus $T^{2}$, but since they are periodic (or almost periodic) the equation 
\eqref{ob} does not apriori constrain the form of $H$ (since the volume form is not determined). 

   Thus, in the case of higher genus, we ignore the relation \eqref{ob} and the related equivalence 
relation \eqref{eq}. Moreover, the conformal group is always compact in this situation, so there 
is no need to divide out by this action. Hence, one works directly with the map
\be \label{Pig}
\Pi: \cH Imm_{+}^{m+1,\a}(\S, \bR^{3}) \to \cC^{m,\a}\times C_{+}^{m-1,\a},
\ee
$$\Pi(F) = ([\g], H).$$
(There is no need to consider the different cases of $\Pi_{0}$ in \eqref{pi01} and $\Pi_{1}$ in \eqref{pi11}). 

   Let $\cM_{c}$ be the Riemann moduli space of constant curvature metrics on the surface $\S$. 
Thus $\g_{c} \in \cM_{c}$ is of constant curvature 0 in case $\S = T^{2}$ and of constant curvature $-1$ 
in case $g > 1$. By the uniformization theorem, any metric $\g$ on $\S$ is of the form 
$\g = \l^{2}\psi^{*}(\g_{c})$, for some diffeomorphism $\psi$. An immersion $F: \S \to \bR^{3}$ 
is conformal if $F^{*}(g_{Eucl})$ is (pointwise) conformal to $\g_{c}$, for some $\g_{c} \in \cM_{c}$. 

  It is then straightforward to verify that all of the results of Section 2 hold for $g > 0$, so that $\Pi$ in 
\eqref{Pig} is a smooth Fredholm map, of index 0. The analog of Proposition 3.1 also holds, although 
the proof requires some further work. The curve of metrics $\g_{t}$ now has the form $\g_{t} = 
\l_{t}^{2}\psi_{t}^{*}(\g_{c(t)})$, where $\g_{c(t)}$ is a curve in $\cM_{c}$. In this situation, the 
relation \eqref{k} for the variation of the induced metric must be replaced by 
$${\tfrac{1}{2}}\k = \d^{*}X^{T} + fA = \f \g + \bar \tau,$$
where $\bar \tau = \l^{2}\tau$ and $\tau$ is tangent to $\cM_{c}$, i.e.~$\tau$ is 
transverse-traceless with respect to the constant curvature metric $\g_{c(t)} \in \cM_{c}$. Following the 
same argument as before, it follows that \eqref{Hest} is modified to 
$$(\int_{\S}H)' = -\int_{\S}H'_{X} - \int_{\S}\<A, \bar \tau\>dV_{\g_{t}},$$
so that 
\be \label{AS}
|(\int_{\S}H)'| \leq K area(\S) + |\int_{\S}\<A, \bar \tau \>|.
\ee
One has 
$$|\int_{\S}\<A, \bar \tau\>| \leq  \int_{\S}|A|^{2} + \int_{\S}|\bar \tau|^{2}.$$
Using \eqref{g} and the Gauss-Bonnet theorem, the first term on the right is bounded by 
$-4\pi\chi(\S) + \int_{\S}H^{2} \leq -4\pi\chi(\S) + H_{0}^{-1}area(\S)$, which is of the same 
form as the first term in \eqref{AS}. Next, since $\g_{t} = \l^{2}\g_{c(t)}$ up to diffeomorphism, 
one has 
$$\int_{\S}|\bar \tau|^{2}dV_{\g_{t}} = \int_{\S}|\tau|^{2}dV_{\g_{c(t)}},$$
where the norm and volume form on the left are with respect to $\g_{t}$ and with respect to 
$\g_{c(t)} \in \cM_{c}$ on the right. However, this term is bounded, since the curve $F_{t}$ has 
bounded speed and length in the target space $\cC^{m,\a}$ and hence in $\cM_{c}$. Thus, both 
the area and the pointwise norm $|\tau|^{2}$ are bounded with respect to $\g_{c(t)}$. It follows 
that \eqref{Hest2} again remains valid in this situation, and the proof is completed as before. 

  The proof of Theorem 3.3 carries over to the higher genus case without change. Thus the map 
$\Pi$ in \eqref{Pig} is smooth, proper and Fredholm, of index 0. 

  However, the computation of the degree $deg \, \Pi$ does not carry over, and it is an open question to 
compute the degree when $g > 0$. (The degree may also depend on the component of the space 
$\cH Imm^{m+1,\a}(\S, \bR^{3})$ if this space is not connected). 

  It would be most natural to compute the degree based on CMC immersed surfaces of higher 
genus in $\bR^{3}$, as done in the case of $S^{2}$. This would require understanding the conformal rigidity 
and infinitesimal conformal rigidity of such CMC immersed surfaces in $\bR^{3}$. 
 
\medskip 

    Next, consider conformal immersions with prescribed mean curvature into the simpy connected 
spaces of constant curvature, i.e.~$\bS^{3}$ and $\bH^{3}$, up to scaling. Theorem 1.1 generalizes 
to this setting with only minor changes. For hyperbolic target $\bH^{3}$, note that $H = 2$ is the 
large-radius limit of the mean curvature of geodesic spheres in $\bH^{3}$, (in place of $H = 0$ 
in $\bR^{3}$). Thus let $C_{2}^{m-1,\a}$ be the space of $C^{m-1,\a}$ functions $H$ on $S^{2}$ 
with $H > 2$ everywhere and let $\cD_{2}^{m-1,\a}$ be the quotient under the equivalence 
relation \eqref{eq} as before. 

\begin{theorem}
The map
\be \label{Pi11}
\Pi_{1}: \cH Imm_{+}^{m+1,\a}(S^{2}, \bS^{3}) \to \cC^{m,\a} \times \cD^{m-1,\a},
\ee
$$\Pi_{1}(F) = ([F^{*}(g_{+1})], [H_{F}]),$$ 
is a smooth and proper Fredholm map of index 0 and
$$deg \, \Pi_{1} = 1.$$
Hence Theorem 1.1 holds with target $\bS^{3}$. 

   Similarly, for hyperbolic space $\bH^{3}$, the map 
\be \label{pi12}
\Pi_{1}: \cH Imm_{+}^{m+1,\a}(S^{2}, \bH^{3}) \to \cC^{m,\a} \times \cD_{2}^{m-1,\a},
\ee
$$\Pi_{1}(F) = ([F^{*}(g_{-1})], [H_{F}]),$$ 
is a smooth and proper Fredholm map of index 0 and
$$deg \, \Pi_{1} = 1.$$
Hence Theorem 1.1 holds with target $\bH^{3}$. 
\end{theorem}

  {\bf Proof:} It is straightforward to verify that all of the discussion and results in Section 2 carry over 
to these target spaces. A simple exercise shows that \eqref{index} remains valid, although it also follows 
from the invariance of the Fredholm index under continuous deformations. 

  The basic method of proof of Propositions 3.1 and 3.2 also carries over, with one difference however. 
Namely, the formula \eqref{HX} for $H'_{X}$ is altered by the presence of curvature to
$$H'_{X} = -\D f - (|A|^{2} + Ric(N,N))f + X^{T}(H).$$
Here $Ric(N,N) = 2\k$, where $\k = \pm 1$ according to whether the target is $\bS^{3}$ or $\bH^{3}$. 
The divergence constraint \eqref{cg} is unaltered, since $Ric(N, X^{T}) = 0$. As before, one then obtains
\be \label{Hk}
(\int_{S^{2}}H)' \leq K area(F_{t}) - 4\k\int_{S^{2}}f.
\ee

 Consider first $\k = 1$, so we are working with the case $\bS^{3}$. One has
$$\int_{S^{2}}f = \int_{S^{2}}\<X, N\>,$$
where $N$ is the outward unit normal. Suppose for the moment the immersions $F_{t}$ extend to 
immersions of a 3-ball $F_{t}: B^{3} \to \bR^{3}$. Let $g_{t} = F_{t}^{*}(g_{+1})$ be the resulting 
curve of constant curvature $+1$ metrics on $B^{3}$. Then
\be \label{fa}
\int_{S^{2}}\<X, N\> = \frac{d}{dt}vol(B^{3}, g_{t}).
\ee
Substituting this in \eqref{Hk} and integrating as before gives
\be \label{vf}
\int_{S^{2}}H \leq K\int_{0}^{t}area(F_{t}) - 4vol(B^{3}, g_{t}) + c \leq  K\int_{0}^{t}area(F_{t}) + c.
\ee
Thus, Proposition 3.1 follows as previously. 

  In general, we may suppose that the initial map $F_{0}$ is an embedding and, by a small perturbation, 
that the normal variation $f$ of $F_{t}$ vanishes only on sets of area zero on $S^{2}$. The embedding 
$F_{0}$ extends to an embedding of the $3$-ball $F_{0}: B_{1} \to \bR^{3}$, with $F_{0}^{*}(g_{+1})$ 
a constant curvature $+1$ metric on $B_{1}$ inducing the metric $\g$ on $\partial B_{1}$. The smooth 
family of mappings $F_{t}: S^{2} \to \bR^{3}$, $0 \leq t \leq T$, then gives a map $F_{T}: B_{T+1} \to \bR^{3}$ 
whose restriction to $S^{2}_{t+1}$ is the immersion $F_{t}$ - so one has a one parameter 
smooth family of immersions of the spheres (a regular homotopy). The pullback $F_{T}^{*}(g_{+1})$ 
is a constant curvature $+1$, possibly singular metric on $B_{T+1}$, but is regular almost 
everywhere. The singular set correponds to the locus where $f = 0$ and so the volume form vanishes. 
Hence \eqref{fa} remains valid, and its integrated version in \eqref{vf} holds for all 
$t$. Thus again Proposition 3.1 follows as before. 

   In the hyperbolic case $\bH^{3}$, the same argument as above gives
$$\int_{S^{2}}H \leq K\int_{0}^{t}area(F_{t}) + 4vol(B^{3}, g_{t}) + c,$$
where the metrics $g_{t}$ are now hyperbolic, i.e.~of constant curvature $-1$. By a well-known 
isoperimetric inequality for hyperbolic metrics, cf.~\cite{Y1} for instance, 
$$vol(B^{3}, g_{t}) \leq \frac{1}{2}area(\partial B^{3}, g_{t}) = \frac{1}{2}area(F_{t}),$$
and hence
$$\int_{S^{2}}H \leq K\int_{0}^{t}area(F_{t}) + 2area(F_{t}) + c.$$
Since $H \geq 2 + H_{0}$ with $H_{0} > 0$, one can absorb the term on the right into the left 
and proceed as before. 

  Thus Proposition 3.1 carries over to both target spaces $\bS^{3}$ and $\bH^{3}$. The proof 
of the other results and Proposition 3.4 carries over to this setting with only very minor changes 
which, as before, completes the proof. 

{\endproof}

  Finally, the discussion of surfaces of higher genus immersed in $\bR^{3}$ carries over without 
further changes to $\bS^{3}$ and $\bH^{3}$. Again, the main remaining question is to compute 
the associated degree.

\bibliographystyle{plain}

\end{document}